\documentclass[12pt]{amsart}

\usepackage{amsmath}
\usepackage {amssymb}
\usepackage{amsfonts}
\usepackage{amscd}
\usepackage{graphicx}
\usepackage{epstopdf}

\def\Z{\mathbb{Z}}
\def\S{\mathbb{S}}

\def\Q{\mathbb{Q}}

\def\bs{\boldsymbol}

\def\a{\alpha}
\def\b{\beta}

\def\g{\gamma}

\def\vempty{\varnothing}
\def\d{\partial}
\def\h{\mathfrak{h}}
\def\<{\langle}
\def\>{\rangle}
\def\V{H_1[k]}
\def\M{H_2[k]}
\def\cl{\overline}

\newtheorem{thm}{Theorem}[section]
\newtheorem{lem}[thm]{Lemma}
\newtheorem{defn}[thm]{Definition}
\newtheorem{prop}[thm]{Proposition}

\newtheorem{claim}[thm]{Claim}
\newtheorem{conj}[thm]{Conjecture}
 
\begin{document}

\title{lens space surgeries \& primitive/Seifert type constructions \ \ }
\author{Michael J. Williams}
\date{\today}
\address{\hskip-\parindent
	Michael J. Williams\\
	Mathematics Department\\
	University of California\\
	Davis, CA 95616\\
	USA}
\email{mikew@math.ucdavis.edu}
\thanks{Research supported in part by NSF VIGRE Grant No. DMS-0135345.}

\begin{abstract}
We show that lens space surgeries on knots in $\S^3$ which arise from the primitive/Seifert type construction also arise from the primitive/primitive construction. This is the first step of a three step program to prove the Berge conjecture for tunnel number one knots.
\end{abstract}	

\maketitle

\section{Introduction}
 
Let $k$ be a knot in the 3-sphere $\S^3$, and let $E(k)$ denote the exterior $E(k)=\cl{\S^3 - k}$. Let $r$ be an isotopy class of simple closed curves on the torus $\d E(k)$. The class $r$ is called a \textbf{slope} of $\d E(k)$. Let $k(r)$ denote the manifold obtained by Dehn surgery on $k$ with slope $r$. As in \cite{ro:knots}, we parametrize slopes of $\d E(k)$ by elements in $\Q \cup \{\frac{1}{0}\}$; we consider $r$ as a ratio $\frac{m}{l}$ where $m,l \in \Z$ with $\gcd(m,l)=1$. In particular, $k(\frac{1}{0})=\S^3$, and if $k_u$ denotes the unknot and $p \ne 0,1$, then $k_u(\frac{p}{q})$ is homeomorphic to the lens space $L(p,q)$. 

There has been an active investigation of which nontrivial knots in $\S^3$ admit lens space surgeries. One of the first papers on the subject is \cite{mo:elem}, in which manifolds obtained by Dehn surgeries on \textit{torus knots} were determined. For a $(p,q)$--torus knot $k$, Moser found that if $m,l \in \Z$ satisfy $|m-lpq|=1$, then $k(\frac{m}{l})$ is homeomorphic to the lens space $L(m,lq^2)$. In particular, $k$ admits integral and non-integral lens space surgeries; the integral surgery slopes are $\frac{m}{l}=\frac{pq \pm 1}{1}$. Moser conjectured that only torus knots admit lens space surgeries. 

A counterexample to Moser's conjecture was found in \cite{br:lens}. There it was shown that if $k$ is the $(-11,2)$-cable on the $(-3,2)$--torus knot, then $k(-23)$ is homeomorphic to the lens space $L(-23,16)$. In \cite{fs:lens}, it was proved that this example generalizes to the collection $\mathcal{C}$ of $(2pq \pm 1,2)$--cables on $(p,q)$--torus knots: the surgery slope must be $\frac{4pq \pm 1}{1}$ and the surgery manifold is homeomorphic to the lens space $L(4pq \pm 1,4q^2)$. In \cite{go:satellite}, it was shown that $\mathcal{C}$ contains all cable knots which admit lens space surgeries. At last in \cite{wu:cyclic}, it was proved that $\mathcal{C}$ actually contains all \textit{satellite knots} which admit lens space surgeries. Also see \cite{bl:lens}. Note that the surgery slope in all of these examples is an integer. The Cyclic Surgery Theorem of \cite{cgls:cyclic} implies that if $k$ is a non-torus knot in $\S^3$ which admits a lens space surgery, then the surgery slope must be an integer and $k$ admits at most two such surgeries. 

Examples of \textit{hyperbolic knots} which admit lens space surgeries were discovered in \cite{fs:lens}. The $(-2,3,7)$ pretzel knot, or ``Fintushel-Stern knot" admits two lens space surgeries. 

By carefully positioning a knot $k$ on a genus 2 Heegaard surface $F$ for $\S^3$, Berge constructed an infinite family of knots in $\S^3$ which admit lens space surgeries \cite{berge:some}. Berge's construction for these lens space surgeries is called the \textit{primitive/primitive} or \textit{$p/p$ construction}. All non-hyperbolic knots which admit lens space surgeries are covered by Berge's constructions, as well as the hyperbolic knots appearing in \cite{fs:lens}. It has been conjectured \cite{ki:probs} that Berge's examples are the only knots in $\S^3$ which admit lens space surgeries. For some recent advances related to this, see \cite{os:lens}, \cite{ni:fibered}, \cite{is:lens}, \cite{dm:lens}, \cite{baker:small}, \cite{saito:lens}, and \cite{bgh:lens}.

A similar construction to Berge's giving knots with Seifert-fibered surgeries was  given by Dean in \cite{de:sfs}. This is the \textit{primitive/Seifert} or \textit{$p/S$ construction}. A slight variation of Dean's construction which also yields knots with Seifert-fibered surgeries is called the \textit{primitive/Seifert-m} or \textit{$p/Sm$  construction}. We will refer to these two constructions as \textit{primitive/Seifert type constructions}. Sometimes the Seifert-fibered space obtained happens to be a lens space, and examples of this phenomenon can be found in \cite{eu:sfs}. We investigate this phenomenon throughout this paper. 

We show that it is always the case that when a lens space surgery arises from a primitive/Seifert type construction, then it also arises from a primitive/primitive construction. 

\begin{thm} \label{thm:ps} 
Let $k$ be a nontrivial knot on a genus 2 Heegaard surface $F$ for $\S^3$. Suppose that $k$ is in $p/S$ position with surface slope $r$. If $k(r)$ is a lens space, then $k$ admits a $p/p$ position with surface slope $r$.
\end{thm}

\begin{thm} \label{thm:psm}
Let $k$ be a nontrivial knot on a genus 2 Heegaard surface $F$ for $\S^3$. Suppose that $k$ is in $p/Sm$ position with surface slope $r$. If $k(r)$ is a lens space, then $k$ admits a $p/p$ position with surface slope $r$. Moreover, $k$ is a torus knot.
\end{thm}

Our approach to Theorem \ref{thm:ps} is to analyze genus 2 Heegaard splittings of Seifert-fibered spaces with two exceptional fibers and the disk as the base space. We rely on the structure theorems in \cite{brz:sfs}. This will enable us to show that the dual knot $k'$ for $k$ in the lens space $k(r)$ is a $(1,1)$ knot in the sense of \cite{doll:bridge}, that is, $k(r)$ has a genus 1 Heegaard splitting $k(r)=V_1 \cup_\Sigma V_2$ such that $k' \cap V_i$ is a trivial arc for each $i=1,2$. This suffices to prove the theorem.

Our approach to Theorem \ref{thm:psm} is to show that $k$ admits a $p/S$ position with the same surface slope. This is done via double branched coverings following the approach in \cite{eu:sfs}. Then we can use techniques from the proof of Theorem \ref{thm:ps} to deduce that $k$ must actually be a torus knot.

This paper is organized as follows. Section 2 reviews the construction of totally orientable Seifert-fibered spaces in terms of Dehn filling on trivial circle bundles. In section 3, we review the surface slope surgery constructions which include the primitive/primitive and primitive/Seifert type constructions. In section 4, we show that the surface slope surgery construction is an effective way to analyze integral surgeries on tunnel number one knots. In section 5, we prove Theorem \ref{thm:ps}. In section 6 we prove Theorem \ref{thm:psm}. 

I would like to thank my advisor Abigail Thompson for her guidance and support. I would also like to thank Kenneth Baker for answering many questions that I had about Berge knots, and for suggesting the claim in the $HD_0$ case of Theorem \ref{thm:ps}. Finally, I would like to thank Cameron Gordon for giving a wonderful and inspiring mini-course on Dehn surgery at Park City, Utah in the summer of 2006, and for helpful conversations.  

\section{Model Seifert Fiberings}

We call an orientable Seifert-fibered space \textbf{totally orientable} if it has an orientable base space. We recall the construction of the \textbf{model Seifert fiberings} following the discussion in \cite{ha:three}. Let $F_0$ be a compact, oriented surface with $m > 0$ boundary circles $c_1,\dots,c_m$. Let $\frac{p_1}{q_1},\dots,\frac{p_k}{q_k} \ (1 \le k \le m)$ be rational numbers in lowest terms. We start with the manifold $M_0=F_0 \times \S^1$ with the natural Seifert fibration where the fibers are $p \times \S^1$ for each $p \in F_0$; we call this the \textbf{trivial Seifert fibration}. Let $T(i)$ correspond to the $i^{th}$ boundary torus $c_i \times \S^1$ of $M_0$. We extend the Seifert fibration of $M_0$ as follows. For each $i$ ($1 \le i \le k$), let $t_i \subset T(i)$ be a an oriented regular fiber so that $t_i \cap c_i=1$ algebraically. This allows us to parametrize slopes on $\d M_0$ using the basis $\{[c_i], [t_i]\}$ for $\mathrm{H}_1(T(i);\Z)$. Now fill each of these $T(i)$ with a solid torus $V_i$ so that the meridian of $V_i$ represents the class $p_i[t_i]+q_i[c_i]$ in $\mathrm{H}_1(T(i);\Z)$. This uniquely determines a Seifert fibering for $V_i$ compatible with $T(i)$. Therefore, we obtain a Seifert-fibered space with base space $F=F_0 \cup_{\d} (k \ \text{disks})$. Let $F(\frac{p_1}{q_1},\dots,\frac{p_k}{q_k})$ denote the resulting Seifert-fibered space. 

Two Seifert-fibered spaces are \textbf{isomorphic} or of the \textbf{same type} if there is a fiber-preserving homeomorphism between them. Note that a product $F \times \S^1$ with the trivial Seifert fibration is of type $F(0)$. In general, we see that $F(\frac{p_1}{q_1},\dots,\frac{p_k}{q_k})$ and $F(\frac{p_1}{q_1},\dots,\frac{p_k}{q_k},0)$ are isomorphic models. We include the following proposition for completeness. 

\begin{prop}{\rm{(Totally orientable case of Proposition 2.1 in \cite{ha:three}.)}} \label{prop:ha} 
Every totally orientable Seifert-fibered space $M$ is isomorphic to some model  $F(\frac{p_1}{q_1},\dots,\frac{p_k}{q_k})$. Two models $F(\frac{p_1}{q_1},\dots,\frac{p_k}{q_k})$ and $F(\frac{p'_1}{q'_1},\dots,\frac{p'_k}{q'_k})$ are isomorphic by an orientation-preserving homeomorphism if and only if, after possibly permuting indices, $\frac{p_i}{q_i} \equiv \frac{p'_i}{q'_i} \ \mathrm{mod} \ 1$ (for each $i$) and, if $\d F=\vempty$, $\sum_i \frac{p_i}{q_i} = \sum_i \frac{p'_i}{q'_i}$. 
\end{prop}

Also see \cite[Theorem 5]{seifert:sfs}. We note that $\S^3$, $\S^1 \times \S^2$, and the lens spaces are the only 3-manifolds which admit Seifert fibrations of type $\S^2(\frac{p_1}{q_1},\frac{p_2}{q_2})$.

\section{Surface Slope Surgery Constructions}

Suppose $F$ is a genus 2 Heegaard surface in $\S^3$. $F$ separates $\S^3$ into genus 2 handlebodies $H_1$ and $H_2$. Suppose that $k$ is a simple closed curve on $F$. Let $N(k)$ denote a tubular neighborhood of $k$; $N(k) \cap F$ is an annulus $A$. Let $E(k)$ denote the knot exterior $\cl{\S^3 - N(k)}$. If we set $F_0=\cl{F - A}$, it is easy to see that $E(k)=H_1 \cup_{F_0} H_2$. The knot $k$ picks up an integral slope $r$ from $F$, called the \textbf{surface slope}. 

We view the $r$-Dehn filling on $E(k)$ as the result of the following procedure. See Figure 1. 

\begin{enumerate}
\item Attach a 2-handle to $H_1$ along $k$ to obtain a $3$--manifold $\V$.
\item Attach a 2-handle to $H_2$ along $k$ to obtain a $3$--manifold $\M$.
\item Identify the boundaries of $\V$ and $\M$ in a way that extends the identification of $H_1$ and $H_2$ along $F_0$.
\end{enumerate}

\begin{figure}[htbp]
  \centering
  \includegraphics[width=3in]{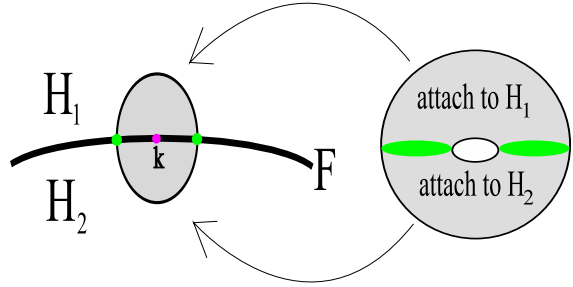}
  \label{surface_slope_surgery}
\caption{Dehn filling on $E(k)$ as a pair of 2-handle attachments to $H_1$ and $H_2$ along $k$.}
\end{figure}

We give some terminology (as in \cite{eu:sfs}) for the position of $k$ with respect to the Heegaard splitting $(H_1,H_2,F)$ in terms of 2-handle attachments.

\begin{defn} Let $k$ be a knot on a genus 2 Heegaard surface $F$ in a splitting $(H_1,H_2,F)$ of $\S^3$. Let $i \in \{1,2\}$.
\begin{itemize}
\item If $H_i[k]$ is a solid torus, then we say that $k$ is \textbf{primitive} on $H_i$. 
\item If $H_i[k]$ admits a Seifert fibration over the disk with at most two exceptional fibers, then we say that $k$ is \textbf{Seifert} on $H_i$.\item If $H_i[k]$ admits a Seifert fibration over the M\"obius band with at most one exceptional fiber, then we say that $k$ is \textbf{Seifert-m} on $H_i$.
\end{itemize}
\end{defn}

We remark that $H_i[k] \cong \S^1 \times D^2$ is equivalent to the property that $k$ meets some meridian disk of $H_i$ in one point. Another equivalent property is that $k$ is part of a free basis of $\pi_1(H_i)$. The way to think about this property is that $k$ is a \textit{core} of $H_i$. 

\begin{defn} Let $k$ be a knot on a genus 2 Heegaard surface $F$ in a splitting $(H_1,H_2,F)$ of $\S^3$. Suppose that $k$ is primitive on $H_1$. Let $r$ be the surface slope.
\begin{itemize}
\item If $k$ is primitive on $H_2$, then we say that $k$ is in \textbf{$\bs{p/p}$ position} on $(H_1,H_2,F)$ and that $k(r)$ arises from the \textbf{primitive/primitive construction}.
\item If $k$ is Seifert on $H_2$, then we say that $k$ is in \textbf{$\bs{p/S}$ position} on $(H_1,H_2,F)$ and that $k(r)$ arises from the \textbf{primitive/Seifert construction}.
\item If $k$ is Seifert-m on $H_2$, then we say that $k$ is in \textbf{$\bs{p/Sm}$ position} on $(H_1,H_2,F)$, and that $k(r)$ arises from the \textbf{primitive/Seifert-m construction}. 
\end{itemize} 
\end{defn}

Berge \cite{berge:some} analyzed the primitive/primitive construction, and his results suggested the following conjecture which appears as \cite[Problem 1.78]{ki:probs}.
\vskip.5cm

\begin{conj}{\rm{(Gordon)}} \label{conj:berge}
If Dehn surgery on a knot $k$ yields a lens space, then $k$ is a Berge knot.
\end{conj}

We remark that Conjecture \ref{conj:berge} asserts only that $k$ admits a primitive/primitive position. This induces a lens space surgery $k(r)$ on $k$. It is possible that $k$ admits another lens space surgery $k(s)$, and it is not clear that $k$ has an alternate primitive/primitive position that realizes this other surgery. For a $(p,q)$--torus knot $k$, both integral lens space surgeries $k(pq \pm 1)$ are easily obtainable by primitive/primitive constructions. As remarked in the introduction, a satellite knot admits at most one lens space surgery, and such a surgery can realized by a primitive/primitive construction. A hyperbolic example is the $(-2,3,7)$ pretzel knot; denote this knot by $k_p$. It is well known that $k_p(18)$ and $k_p(19)$ are lens spaces. Berge remarks in \cite{berge:solid} that this knot appears to be embeddable in an unknotted solid torus $V \subset \S^3$ so that two surgeries on $k_p \subset V$ will yield a solid torus; it is an exercise to show that this is actually the case. Then by Gabai's Theorem on knots in solid tori \cite[Theorem 1.1]{ga:solid}, we conclude that $k_p$ and its dual knots in $k_p(18)$ and $k_p(19)$ are indeed (1,1) knots.  

We actually expect that any lens space surgery $k(r)$ arises from a primitive/primitive construction. In the special cases of Theorem \ref{thm:ps} and \ref{thm:psm}, this stronger condition is satisfied.  

The primitive/Seifert construction is studied in the papers \cite{de:sfs}, \cite{eu:sfs}, and \cite{mm:sfs}. In \cite{eu:sfs} and \cite{mm:sfs}, the primitive/Seifert-m construction is also studied. 

\section{Integral Surgery on Tunnel Number One Knots}

Let $k_0$ be a tunnel number one knot in $\S^3$, and let $E(k_0)$ denote the exterior of $k_0$ in $\S^3$. We can realize any integral surgery on $k_0$ as a surface slope surgery on a genus 2 Heegaard surface as follows. Let $t$ be a tunnel for $k_0$. Let $H_1$ be a tubular neighborhood of $k_0 \cup t$ and set $H_2=\cl{\S^3 - H_1}$. If we also set $F=\d H_1=\d H_2$, then $(H_1,H_2,F)$ is a genus 2 Heegaard splitting of $\S^3$. Note that $k_0$ is a core of $H_1$. 

Let $r \in \Z$ be given. Now push $k_0$ through $H_1$  to the boundary $F$ in order to obtain a copy $k$ of $k_0$ that has surface slope $r$. Note that $k$ is primitive on $H_1$. See Figure 2.

\begin{figure}[htbp]
  \begin{center}
  \includegraphics[width=3in]{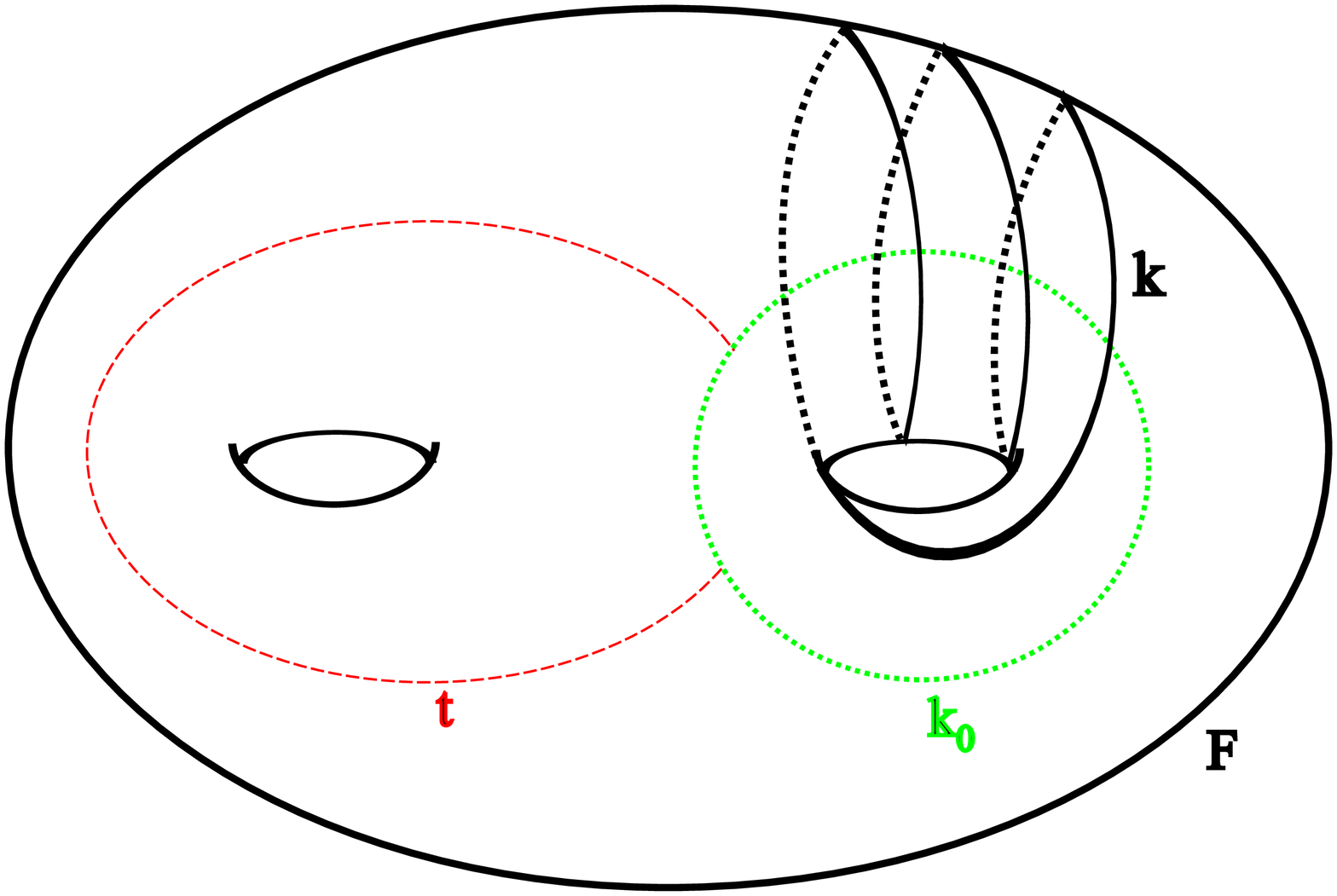}
  \label{push_off}
\caption{Pushing off $k_0$ to a copy $k \subset F$ with surface slope $r$.}
  \end{center}
\end{figure}  

Now suppose that $k_0$ is hyperbolic and $k_0(\a)$ is a lens space $M$ for some $\a \in \Q$. The Cyclic Surgery Theorem of \cite{cgls:cyclic} implies that $\a \in \Z$ since $k_0$ is not a torus knot. Therefore we can realize the surgery slope by the above construction with $r=\a$. Note that $|r| \ge 2$. If $\M$ is a solid torus, then $k$ is in $p/p$ position. 

Recall that a 3-manifold $Y$ with $\d Y \ne \vempty$ is \textbf{$\bs{\d}$--reducible} if $\d Y$ compresses in $Y$; otherwise,  $Y$ is \textbf{$\bs{\d}$--irreducible}. Also, if $\g_1$ and $\g_2$ are simple closed curves on a torus, we let $\Delta(\g_1,\g_2)$ denote their minimal geometric intersection number. 

The following lemma tells us some important properties of  $\M$.

\begin{lem} \label{lem:haken}
Let $k_0$ be a nontrivial knot in $\S^3$ with unknotting tunnel $t$. Let $H_1$ be a regular neighborhood of $k_0 \cup t$ and let $H_2$ be $\cl{\S^3 - H_1}$. Push $k_0$ through $H_1$ to a copy $k \subset \d H_1$. Let $r \in \Z$ be the resulting surface slope. Assume $k(r)$ is irreducible. Then $\M$ is irreducible. Furthermore, $\M$ is $\d$--reducible if and only if $\M \cong \S^1 \times D^2$.
\end{lem} 

\begin{proof}
There are two cases according to whether or not $\d H_2 - k$ is compressible in $H_2$. Suppose that $\d H_2 - k$ is incompressible in $H_2$. Since $H_2$ is irreducible and has compressible boundary, the work of \cite{jaco:handle} (see also \cite{prz:incompressible} and \cite{cg:reducing}) asserts that $\M$ is irreducible and has incompressible boundary. 

Now suppose that $\d H_2 - k$ has a compressing disk $D$. Cut $H_2$ along $D$ to obtain a 3-manifold $N$. We can assume that $D$ is nonseparating so that $N$ is a solid torus. Let $D'$ be a meridian disk of $N$. Since $\cl{H_2 - (D \cup D')}$ is a 3-ball and $k$ is nontrivial, $\Delta(k,\d D') \ge 1$. Let $u$ denote the unknot. It is not difficult to see that $\M = (\S^1 \times D^2) \# k_u(\frac{a}{b})$ where $a=\Delta(k,\d D')$ and $b$ is some integer coprime to $a$. In particular, $\M \cong \S^1 \times D^2$ if and only if $a=1$. Now suppose that $a > 1$. Hence $k_u(\frac{a}{b})=L(a,b)$. Let $\mu$ be a meridian of the tunnel $t$ for $k$. Then $\mu$ is meridian for the solid torus $\V$. Now, $k(r)=\V \cup_\d \M=(\S^1 \times D^2)  \cup_\d [(\S^1 \times D^2) \# L(a,b)]$. Unless $\Delta(\mu,\d D)=1$, $k(r)$ will be the (nontrivial) connect-sum of lens spaces, contradicting irreducibility of $k(r)$. This implies that the genus 2 Heegaard splitting $(\cl{H_1 - k_0}, H_2, F)$ for $E(k_0)$ is stabilized, contradicting that $k_0$ is nontrivial. This shows that $\M$ is $\d$--reducible if and only if $\M \cong \S^1 \times D^2$. 
\end{proof}

By Lemma \ref{lem:haken}, we can assume that $\M$ is irreducible. Then by Thurston's Hyperbolization Theorem for Haken 3-manifolds \cite{th:three}, the three possibilities for the topology of $\M$ are: Seifert-fibered, toroidal, or hyperbolic (more precisely, $int(\M)$ admits a hyperbolic metric). This paper examines the case in which $\M$ is Seifert-fibered. The case in which $\M$ is toroidal is investigated in \cite{wi:pt}. 

It is well-known that any tunnel number one knot in $\S^3$ is strongly invertible. The problem of obtaining Seifert-fibered spaces via Dehn surgery on hyperbolic, strongly invertible knots has been studied by Eudave-Mu\~noz in \cite{eu:sfs}. Each of the knots studied in \cite{eu:sfs} can be made to lie as a non-separating curve on a genus 2 Heegaard surface in $\S^3$; the resulting surface slope surgeries are studied. In the case that such a surgery is atoroidal, the following proposition asserts that it suffices to study tunnel number one knots. 

\begin{prop}
Let $k$ be a hyperbolic knot which lies on a genus 2 Heegaard surface $F$ for $\S^3$, and let $r$ denote the resulting surface slope. If $k(r)$ is atoroidal, then $k$ has tunnel number one. Moreover, $k$ is primitive on one of the handlebodies bounded by $F$.  
\end{prop}

\begin{proof}
Let $H_1$ and $H_2$ be genus 2 handlebodies for the induced Heegaard splitting for $\S^3$ given by $F$.  We write this as $\S^3=H_1 \cup_F H_2$. Let $T=\d H_1[k]=\d H_2[k]$. Then $k(r)=H_1[k] \cup_T H_2[k]$. 

We claim that $F-k$ compresses in $H_i$ for some $i$. Suppose, on the contrary, that $F-k$ were incompressible in both $H_1$ and $H_2$. Then each $H_i[k]$ would be $\d$--irreducible by the work of \cite{jaco:handle} (see also \cite{prz:incompressible} and \cite{cg:reducing}). Thus $T$ would be an incompressible torus in $k(r)$, contradicting our assumption that $k(r)$ is atoroidal.

Now suppose that $D$ is an essential disk in, say, $H_1$ which misses $k$. Then the 3-manifold obtained by cutting $H_1$ along $D$ contains a solid torus component $N$ with $k \subset \d N$. Let $D'$ be a meridian of $N$. If $\Delta(k,\d D')>1$, then $k$ is a torus knot or a cable knot,  a contradiction. If $\Delta(k,\d D')=0$, then $k$ lies on a 2-sphere, a contradiction. Thus $\Delta(k,\d D')=1$, exhibiting $k$ as primitive on $H_1$. Consequently, $k$ has tunnel number one.
\end{proof}

\section{Primitive/Seifert Constructions which Also Arise from Primitive/primitive Constructions}

We now prove Theorem \ref{thm:ps}, that is, if $k(r)$ is a lens space which arises from Dean's primitive/Seifert construction, then $k(r)$ also arises from Berge's primitive/primitive construction. Examples of knots in $\S^3$ which have $p/p$ and $p/S$ positions yielding the same lens space surgery can be found in \cite{eu:sfs}.

\begin{defn}
Let $M$ be a $\S^3$, $\S^1 \times \S^2$ or a lens space. A knot $k \subset M$ is a \textbf{$\bs{(1,1)}$ knot} if $M$ has a genus 1 Heegaard splitting $M=V_1 \cup_\Sigma V_2$ such that $k \cap V_i$ is a trivial arc for each $i=1,2$. 
\end{defn}

To prove the theorem, we will show that the dual knot $k'$ is a $(1,1)$ knot in $k(r)$. The following lemma asserts that this is sufficient. 

\begin{lem} \label{lem:1-bridge} Let $M'$ be a lens space. Suppose that $k' \subset M'$ is a $(1,1)$ knot. If $k'$ admits an $\S^3$ surgery, then the dual knot $k \subset \S^3$ is a Berge knot.
\end{lem}

\begin{proof}
A proof for this lemma can be found contained within the proof of Theorem 2 of \cite{berge:some}. For the reader's convenience, here is a review of the proof. Let $(V_1,V_2,\Sigma)$ be a genus 1 Heegaard splitting for $M$ so that $k'$ is a $(1,1)$ knot. Thus $k'$ meets $\Sigma$ in two points, and the arcs $\a=k \cap V_1$ and $\b=k \cap V_2$ are unknotted. 

Let $N$ be a small regular neighborhood of $\a$ in $V_1$. Drill out $N$ from $V_1$ to obtain a genus 2 handlebody $H=\cl{V_1 - N}$. Since $N$ is unknotted in $V_1$, there is a \textit{trivializing disk} $D$ for $N$ in $V_1$; that is, $\d D$ lies in  $\d H$, and the arc $\a$ is parallel in $N$ to $N \cap D$. Note that the disk $D$ is a compressing disk for $H$. Let $A$ be the annulus $\cl{\d N - \Sigma}$. Let $m$ be a core of $A$; note that $m$ meets $D$ in a single point. This makes $m$ primitive on $H$. Let $D'$ be a compressing disk for $V_2$ which misses $\b$. Let $E(k')$ denote the exterior of $k'$ in $M'$. Then $E(k')$ is homeomorphic to $H \cup_{\d D'} (\text{2-handle})$. Therefore, $\d H$ is a genus 2 Heegaard surface for any Dehn filling on $E(k')$.

We may as well assume that $E(k')$ is not Seifert-fibered. Let $\g$ be a slope on $\d E(k')$ which yields an $\S^3$ surgery. By the Cyclic Surgery Theorem of \cite{cgls:cyclic}, we must have $\Delta(m,\g)=1$. Let $k$ denote the dual knot of the $\S^3$ surgery on $k'$. Hence $m=k$. The $\S^3$ surgery on $k'$ with slope $\g$ corresponds to the handle decomposition $$\S^3=(H \cup_{\d D'} (\text{2-handle})) \cup_{\g} (\text{2-handle}) \cup (\text{3-handle}).$$ Let $H'$ denote $\cl{\S^3 - H}$. Thus $(H,H',\d H)$ is a genus 2 Heegaard splitting for $\S^3$. Since $\Delta(m,\g)=1$, $k$ is also primitive on $H'$. Therefore, $k$ is $p/p$ in the genus 2 Heegaard splitting $(H,H',\d H)$.
\end{proof}

\begin{proof}[Proof of Theorem \ref{thm:ps}]
By the preceeding lemma, it is enough to show that the dual knot of our surgery is a $(1,1)$ knot in $k(r)$. By assumption, $\M$ admits a Seifert fibration of type $D^2(\frac{p}{q},\frac{a}{b})$. If $\M$ is a solid torus, then there is nothing to prove; so we will assume that $\mathrm{min}\{|q|,|b|\} > 1$.

Genus 2 Heegaard splittings of Seifert-fibered spaces of type $D^2(\frac{p}{q},\frac{a}{b})$ (generalized torus knot exteriors) were classified in \cite{brz:sfs}: any splitting is isotopic to one of three types, called $HD_0$, $HD_S$ and $HD_T$. We will apply their results in the next few subsections. Let $\tau$ denote the co-core of the 2-handle to be attached to $H_2$ to form $\M$. Note that $\tau=k' \cap \M$. Fix a Seifert fibration for $\M$ of type $D^2(\frac{p}{q},\frac{a}{b})$. Let $D$ denote the base disk and let $\pi:\M \rightarrow D$ denote the projection map. We denote the images of the exceptional fibers of multiplicities $|q|$ and $|b|$ by $S$ and $T$ respectively.

\begin{figure}[htbp]
  \centering
  \includegraphics[width=4in]{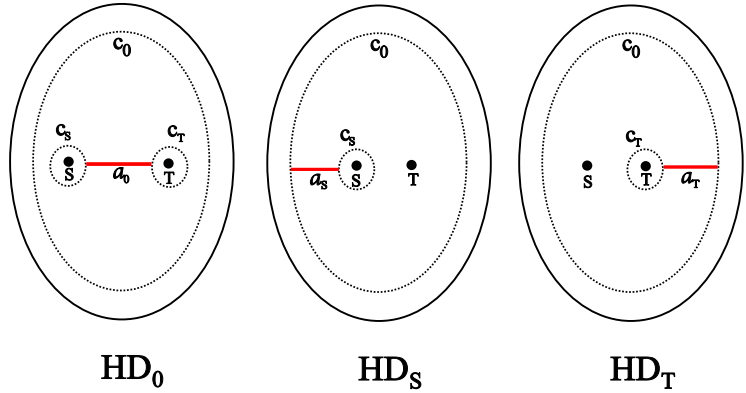}
  \label{vertical_splittings}
\caption{Base space picture of the genus 2 splittings of $\M$.}
\end{figure}  

As described in \cite{brz:sfs}, let $c_0, c_S, c_T$ be disjoint circles in $Int(D)$ so that $\pi^{-1}(c_0)$ is a boundary parallel torus and $\pi^{-1}(c_S)$ (respectively $\pi^{-1}(c_T)$) is a torus bounding a tube about the exceptional fiber $\pi^{-1}(S)$ (respectively $\pi^{-1}(T)$). Let $a_0, a_S, a_T$ be arcs in $Int(D)$ so that $a_0$ connects $c_S$ to $c_T$, $a_S$ connects $c_0$ to $c_S$, and $a_T$ connects $c_0$ to $c_T$. We will now review (as described in \cite{brz:sfs}) the three genus 2 Heegaard splittings $HD_0$, $HD_S$ and $HD_T$ of $\M$ (refer to Figure 3), and show that in each case, the dual knot $k'$ is a $(1,1)$ knot.

\begin{figure}
  \centering
  \includegraphics[width=4in]{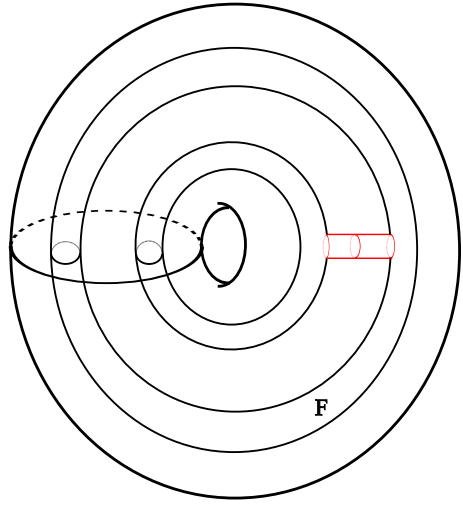}
  \label{HD0_splitting}
\caption{The Heegaard surface $F$ in the splitting $HD_0$.}
\end{figure}
  
\subsection{The Heegaard splitting $HD_0$}
The genus 2 handlebody $H_2$ is obtained by joining the solid tori bounded by $\pi^{-1}(c_S)$ and $\pi^{-1}(c_T)$ by a tube about an essential arc $Z$ in the fibered annulus $\pi^{-1}(a_0)$. See Figure 4. Note that the complement of this handlebody is a genus 2 compression body obtained by attaching a 1-handle $\h$ to $\S^1 \times \S^1 \times [0,1]$. This 1-handle $\h$ will be thought of as a 2-handle when we view it as being attached to $H_2$. The 2-handle $\h$ is a regular neighborhood of the fibered annulus $\pi^{-1}(a_0)$ cut along a neighborhood of $Z$. The co-core of $\h$ is $k' \cap \M$. Note that the attaching curve for this 2-handle intersects one meridian disk of $H_2$ in $|q|$ points and another meridian disk in $|b|$ points.
 
\begin{figure}[htbp]
  \centering
  \includegraphics[width=5in]{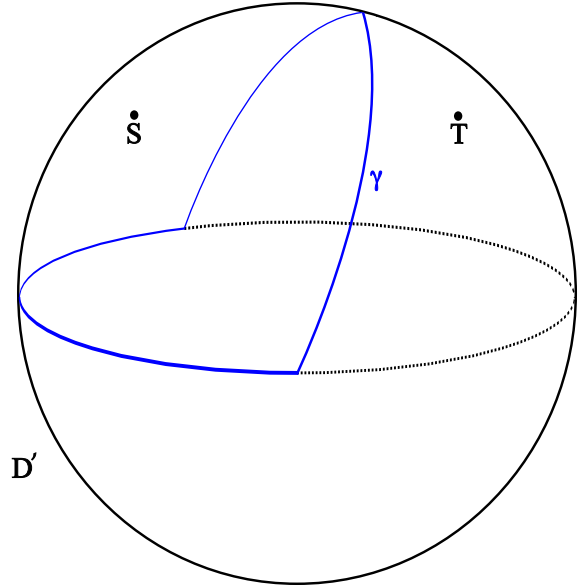}
  \label{choose_gamma}
\caption{Choosing the curve $\g$ on $D'$: The lift of $\g$ in $k(r)$ is a Heegaard torus.}
\end{figure}

Let $M$ denote $k(r)$. Let $D' \cong \S^2$ denote the base surface for the extended Seifert fibration on $M$; let $\pi:M \rightarrow D'$ denote the natural projection obtained by extending the map $\pi:\M \rightarrow D$. 

Let $\g$ be a circle on $D'$ consisting of a proper arc in $D$ which separates $S$ and $T$, and an arc in $\d D$ which is on the side of $S$. See Figure 5. We also assume that $\g$ meets $a_0$ transversely in a single point. Now, $\Sigma=\pi^{-1}(\g)$ is a vertical Heegaard torus for $M$; $M=V_1 \cup_{\Sigma} V_2$, where $V_1$ and $V_2$ are solid tori bounded by $\Sigma$. Label the $V_i$ so  that $S \in \pi(V_1)$ and $T \in \pi(V_2)$. We see (by design) that, in $M$, $\V \cap \Sigma$ in an annulus. We also see that $\h$ intersects $\Sigma$ in a disk containing the co-core of $\h$. Recall that the co-core of $\h$ is $k' \cap \M$. Hence $k' \cap \M$ is a trivial arc in, say, $V_1$. In the other solid torus $V_2$, we can consider $\V$ as part of $V_2$. The arc $k' \cap \V$ is trivial in $\V$, and it does not intersect some meridian disk of $\V$. That is, it does not wind longitudinally around $\V$. Thus $k' \cap V_2$ is a trivial arc in $V_2$ as well as in $\V$. Therefore, $k'$ is a $(1,1)$ knot. Moreover, we have the following. 

\begin{claim} \label{claim:torus1}
$k$ is a torus knot.
\end{claim}

\textit{Proof of claim.} In the step of the proof where we extend the Seifert fibering of $\M$ to one on $M$, we compatibly fiber the solid torus $\V$. Let $D^*$ be a meridian disk of $\V$ containing the arc $k' \cap \V$. We isotop the Seifert fibering on $\V$ so that the every regular fiber intersects $D^*$ exactly once. Hence $A$ meets $D^*$ in a single arc. Let $A^*$ be the complementary annulus $\cl{\Sigma - A}$. Then we may push $k' \cap \V$ through $D^*$ to lie on $A^*$, all without moving $k' \cap \M$. This makes $k'$ lie on $\Sigma$. Therefore, $k'$ is a torus knot in $M$. Hence $k$ is a torus knot in $\S^3$.    

\subsection{The Heegaard splitting $HD_S$}

\begin{figure}
  \centering
  \includegraphics[width=4in]{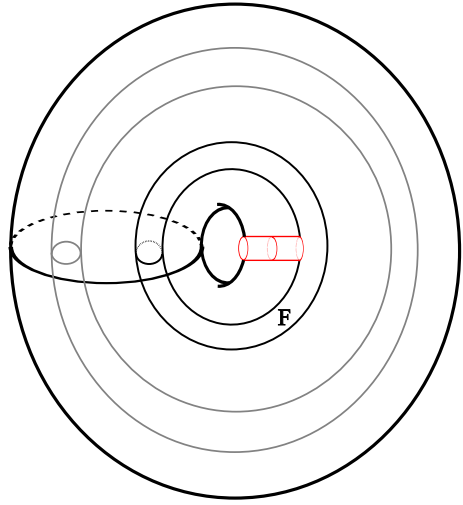}
  \label{HDS_splitting}
\caption{The Heegaard surface $F$ in the splitting $HD_S$.}
\end{figure}

The Heegaard surface $F$ is obtained by connecting the tori $\pi^{-1}(c_S)$ and $\pi^{-1}(c_0)$ by $\d U$ where $U$ is a tube about an essential arc in the fibered annulus $\pi^{-1}(a_S)$. See Figure 6. The surface $F$ bounds the handlebody $H_2$ which has a meridian disk $D_1=\cl{\pi^{-1}(a_S) - U}$. To find another meridian disk $D_2$, cut $H_2$ along $D_1$ to obtain a tube about the critical fiber $\pi^{-1}(T)$, then let $D_2$ be a meridian of this solid torus. The cores of $H_2$ consist of a section over $c_S$ and a core of the solid torus bounded by $\pi^{-1}(T)$. The 2-handle $\h$ attached to $H_2$ in $\M$ is just a product neighborhood of a meridian disk of a solid torus neighborhood of $\pi^{-1}(S)$. See \cite{brz:sfs} for more details.  
  
We see that $k'$ is a $(1,1)$ knot as follows. The arc $\tau=k' \cap \M$ is the co-core of the 2-handle $\h$. Thus, it is an arc of the exceptional fiber $\pi^{-1}(S)$. Now fill $\M$ with $\V$ at $\d \M$. Consider the genus 1 Heegaard splitting $V_1 \cup_{\Sigma} V_2$ of $M$ where $\Sigma=\pi^{-1}(c_S)$ and $S \in \pi(V_1)$. Now, isotope (expand) $V_1$ so that $\pi(\d V_1)=\g$ (see Figure 5). We see that $\Sigma$ meets $\V$ in an annulus. The arc $k' \cap V_1$ is trivial in $V_1$ since $k' \cap V_1$ is contained in a core of $V_1$.  As shown before, $k' \cap V_2$ is trivial. Therefore $k'$ is a (1,1) knot. 

\subsection{The Heegaard splitting $HD_T$}
This is similar to the $HD_S$ case; just replace $S$ with $T$ throughout the proof. \\

This completes the proof of Theorem \ref{thm:ps}.
\end{proof}

\section{Primitive/Seifert-m Constructions which Also Arise from Primitive/primitive Constructions}

We now prove Theorem \ref{thm:psm}, that is, if $k(r)$ is a lens space which arises from the primitive/Seifert-m construction, $k$ is a torus knot. 

\begin{proof}[Proof of Theorem \ref{thm:psm}]
We will show that $k$ admits a $p/S$ position on $F$. Then we can apply techniques from the proof of Theorem \ref{thm:ps}. We may as well assume that $\M$ is not the twisted $I$-bundle over the Klein bottle by Theorem \ref{thm:ps}; this 3-manifold admits a Seifert fibration of type $D^2(\frac{1}{2},\frac{1}{2})$.

We can now set-up our situation as in \cite[Section 3]{eu:sfs}. We consider the genus 2 Heegaard splitting $\S^3=H_1 \cup_F H_2$. If $k$ sits on $F$ with surface slope $r$, then $k(r)=H_1[k] \cup_{\d} H_2[k]$. There is an involution $H_1 \! \rightarrow \! H_1$ with fixed point set consisting of three trivial arcs, two of which have an endpoint lying on $k$.  The involution can be extended to $H_2$ showing that $k$ is a strongly invertible knot. We see that $H_1$ and $H_2$ cover trivial 3-string tangles. Since $k \subset \S^3$, the quotient of $E(k)$ by the strong inversion on $\S^3$ is a 2-string tangle $(B,t)$ which can be summed with a trivial 2-string tangle to form the unknot. The involution on $H_i$ can be extended to $H_i[k]$ for each $i=1,2$. In each case, the fixed point set consists of two arcs and possibly a simple closed curve. 

Since $\V \cong \S^1 \times D^2$, any involution on a solid torus $V$ is fiber preserving \cite[Lemma 6]{to:sfs}; thus it must cover a trivial 2-string tangle. 

Now $\M$ is a Seifert-fibered space over the M\"obius band and the involution preserves fibers on $\d \M$. The main theorem of \cite{to:sfs} implies that the involution is fiber-preserving on all of $\M$. Thus $\M$ covers a Montesinos-m tangle of length one. See Figure 7. 

\begin{figure}
  \centering
  \includegraphics[width=2.5in]{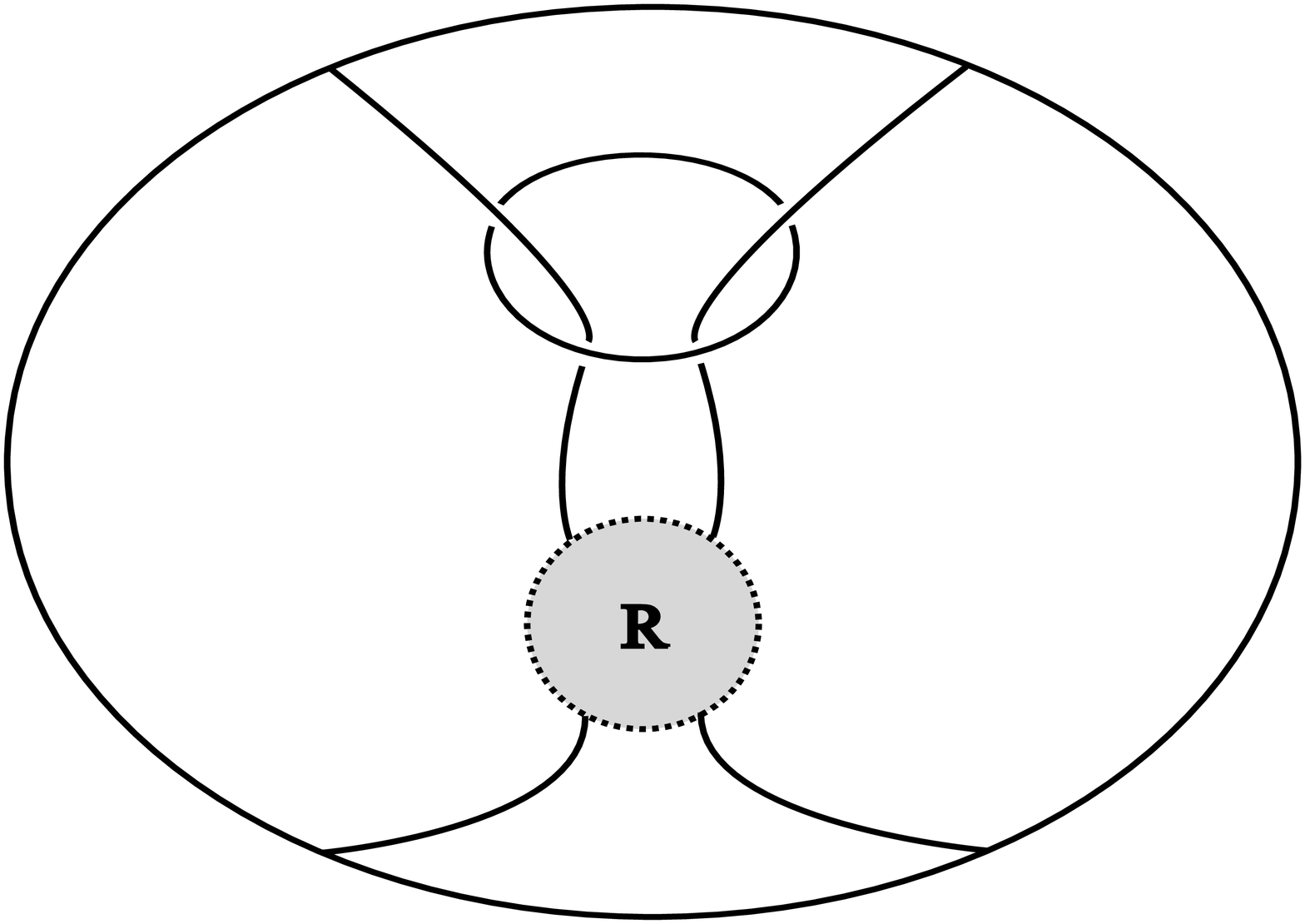}
  \label{Mm}
\caption{A Montesinos-m tangle of length one. The label $\mathrm{R}$ indicates a trivial tangle.}
\end{figure}

In the double branched cover $\pi: k(r) \rightarrow \S^3$, the branch set downstairs is a 2-bridge link $L$. Set $S=\pi(F \cap E(k))$. The surface $S$ is a properly embedded disk in $B$ which intersects $t$ transversally in four points, with none on $\d t$. We sum $(B,t)$ with a trivial tangle $(B',t')$ to form $(\S^3,L)$. Let $s$ denote $\d S$. The disk $S$ meets $B$ in $s$, and the circle $s$ bounds a disk $D' \subset B'$ which separates the strings of $t'$. See Figure 8.

The circle $s$ separates $\d B$ into two disks $D_1$ and $D_2$ where $D_i=\pi(H_i \cap \d E(k))$. The 2-spheres $S \cup D_1$ and $S \cup D_2$ bound 3-balls $B_1$ and $B_2$, respectively, with $B_1 \cap B_2 = S$. For each $i=1,2$, the 3-ball $B_i$ determines a 3-string tangle $(B_i,t_i)$ where $t_i=B_i \cap t$. 

\begin{figure}[htbp]
  \centering
  \includegraphics[width=5in]{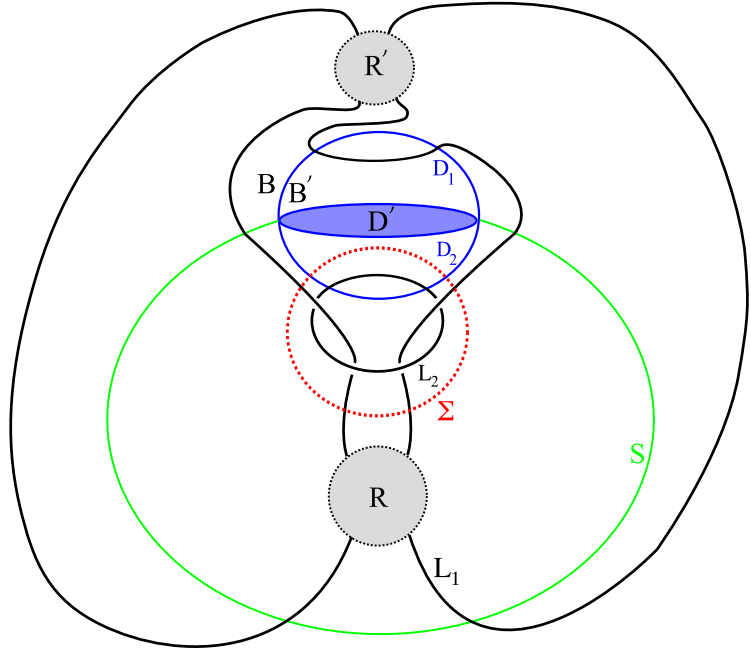}
  \label{tangles}
\caption{A schematic picture of the tangles in the downstairs portion of the double branched cover. The tangles labeled $\mathrm{R}$ and $\mathrm{R}'$ are trivial tangles.}
\end{figure}

We set $\hat{S}=S \cup D'$. Then $\hat{S}$ determines a pair of 2-string tangles $(\hat{B}_1,\hat{t}_1)$ and $(\hat{B}_2,\hat{t}_2)$ where $(\hat{B}_1,\hat{t}_1)$ is a trivial 2-string tangle, and $(\hat{B}_2,\hat{t}_2)$ is a Montesinos-m tangle of length one. 

Let $L_2$ be the component of the 2-bridge link $L$ such that $L_2 \subset \hat{B}_2$ and let $L_1$ be the other component. Note that $L_2$ must consist of an arc from the tangle $(B_2,t_2)$, and an arc from the tangle $(B',t')$. Set $\gamma=L_2 \cap B'=L_2 \cap t'$. 

We assert that there is a properly embedded disk $D$ in $\hat{B}_2$ such that 
\begin{itemize}
\item $\d D = L_2$
\item $D$ meets $L_1$  in two points
\item $D$ meets $B'$ in a disk cut off by $\gamma$
\end{itemize}
Such a disk $D$ exists because we can push $\gamma$ through $\d B'$ so that $\cl{\hat{B}_2 - B'}$ is a Montesinos-m tangle. We now thicken $D$ to a 3-ball with boundary $\Sigma$. See Figure 8. Note that $\Sigma \cap B'$ is a single disk disjoint from $\gamma$. We also note that $\Sigma \cap B'$ is isotopic to $\hat{S} \cap B'$ in $B' - t'$. Without loss of generality, we may assume that $\Sigma \cap D'=\vempty$. 

The sphere $\Sigma$ bounds 3-balls $A_1$ and $A_2$ where $L_2 \subset A_2$. For each $i=1,2$, let $s_i$ denote $A_i \cap L$. Note that $(A_2,s_2)$ is a Montesinos tangle whose lift $\pi^{-1}(A_2)$ is the twisted $I$-bundle over the Klein bottle. This manifold is $\d$--irreducible.

We now show that $(A_1,s_1)$ is a trivial tangle. Note that the lift $\pi^{-1}(\Sigma)$ is a separating compressible torus in the lens space $k(r)$. Since $\pi^{-1}(A_2)$ is $\d$--irreducible, $\pi^{-1}(A_1)$ must be $\d$--reducible. But $(A_1,s_1)$ can be tangle summed with $(A_2,A_2 \cap L_1)$ to form $(\S^3,L_1)$, so $\pi^{-1}(A_1)$ is the exterior of some knot in $\S^3$. Since $\pi^{-1}(A_1)$ is $\d$-reducible, it must be the exterior of the unknot. Hence $\pi^{-1}(A_1)$ is a solid torus. It follows that $(A_1,s_1)$ is a trivial tangle.  

Now push $\Sigma$ so that $\Sigma \cap B'=D'$. Set $\Sigma_0=\cl{\Sigma - B'}$. We see that $\Sigma=\Sigma_0 \cup D'$ where $D'=\hat{S} \cap B'$. For each $i=1,2$, the replacement of $D'$ in $\Sigma_0 \cup D'$ by $D_i$ in $\Sigma_0 \cup D_i$ is accomplished by an isotopy in $B'$. This yields a new $\d$-parallel arc to form the tangle $(A_i',s_i')$ corresponding to $(A_i,s_i)$. This shows that $(A_i',s_i')$ is a trivial 3-string tangle for each $i=1,2$. Thus $\pi^{-1}(A_i')$ is  a genus 2 handlebody in $\S^3$ whose boundary contains $k$. The surface slope is still $r$ because $\Sigma \cap B'=\hat{S} \cap B'=D'$. This shows that $k$ is in $p/S$ position on a genus 2 Heegaard surface for $\S^3$ with surface slope $r$. Moreover, we have the following.

\begin{claim} \label{claim:torus2}
$k$ is a torus knot.
\end{claim}

\textit{Proof of claim.} Let $H_i'=\pi^{-1}(A_i')$ for $i=1,2$. We have just shown that $k$ is primitive on $H_1'$ and Seifert on $H_2'$ and the surface slope is still $r$. We have also shown that $H_2'[k]$ is actually a twisted $I$-bundle over the Klein bottle. This 3-manifold is a Seifert-fibered space of type $D^2(\frac{p_1}{2},\frac{p_2}{2})$ for some odd integers $p_1,p_2$. Then it immediately follows from Corollary 5.8 in \cite{brz:sfs} that there is only one genus 2 Heegaard splitting of $H_2'[k]$ up to isotopy. In particular, we obtain $k' \cap H_2'[k]$ from the $HD_0$ splitting for $H_2'[k]$. As in the proof of Claim \ref{claim:torus1}, it follows that $k$ is a torus knot. This completes the proof of the claim.
\end{proof}

There is an alternate indirect way to prove Claim \ref{claim:torus2}.  By Theorem \ref{thm:ps}, $k(r)$ arises from the double primitive construction. Since $k(r)$ also arises from the $p/Sm$ construction, we have  $k(r) \cong L(4n,2n-1)$ for some $n \ne 0$. A recent result of Ichihara and Saito \cite[Theorem 1.1]{is:lens} can now be applied to conclude that $k$ must be either of the $(\pm5,3)$--, or, the $(\pm7,3)$--torus knots.

\end{document}